\documentclass[12pt]{article}
\date{May 2000}
\usepackage{amssymb}
\input epsf.sty

\def\G{\Gamma}
\def\Om{\Omega}

\def\a{\alpha}
\def\g{\gamma}
\def\om{\omega}

\def\vp{\phi}

\def\bu{\bullet}
\def\ify{\infty}

\def\longra{\longrightarrow}
\def\mpo{\mapsto}
\def\ot{\otimes}
\def\ov{\overline}
\def\ra{\rightarrow}
\def\sbs{\subset}
\def\ts{\times}
\def\wdg{\wedge}
\def\wt{\widetilde}

\font\tenbb=msbm10
\font\sevenbb=msbm7
\font\fivebb=msbm5
\newfam\bbfam
\textfont\bbfam=\tenbb \scriptfont\bbfam=\sevenbb
\scriptscriptfont\bbfam=\fivebb
\def\bb{\fam\bbfam}

\def\sevafig#1#2{\centerline{
 \epsfxsize=#2\epsfbox{#1}}}

\def\Cb{{\bb C}}
\def\Rb{{\bb R}}

\def\Dc{{\cal D}}
\def\Fc{{\cal F}}
\def\Gc{{\cal G}}
\def\Hc{{\cal H}}
\def\Mc{{\cal M}}

\def\Uc{{\cal U}}

\def\build#1_#2^#3{\mathrel{
\mathop{\kern 0pt#1}\limits_{#2}^{#3}}}

\title{Deformation quantization with traces}

\author{Giovanni Felder and Boris Shoikhet}

\begin{document}

\maketitle

\begin{abstract}
In the present paper we prove a statement closely related 
to the cyclic formality conjecture \cite{Sh}. In 
particular, we prove that for a constant volume form $\Om$ and a 
Poisson bivector field $\pi$ on $\Rb^d$ such that ${\rm 
div}_{\Om} \, \pi = 0$, the Kontsevich star-product \cite{K} 
with the harmonic angle function is cyclic, i.e. 
$\int_{\Rb^d} \, (f * g) \cdot h \cdot \Om = \int_{\Rb^d} \, (g * 
h) \cdot f \cdot \Om$ for any three functions $f,g,h$ on 
$\Rb^d$ (for which the integrals make sense). We also prove 
a globalization of this theorem in the case of arbitrary 
Poisson manifolds and an arbitrary volume form,
and prove a generalization of the 
Connes-Flato-Sternheimer conjecture \cite{CFS} on closed 
star-products in the Poisson case.
\end{abstract}

\section{Cyclic formality conjecture}\label{sec1}

We work with the algebra $A = C^{\ify} (M)$ of smooth 
functions on a smooth manifold $M$. One associates to the 
algebra $A$ two differential graded Lie algebras: the Lie 
algebra $T_{\rm poly}^{\bu} (M)$ of smooth polyvector fields 
on the manifold $M$ (with zero differential and the 
Schouten-Nijenhuis bracket), and the polydifferential part 
$\Dc_{\rm poly}^{\bu} (M)$ of the cohomological Hochschild 
complex of the algebra $A$, equipped with the Gerstenhaber 
bracket (see \cite{K} for the definitions). We consider 
$T_{\rm poly}^{\bu} (M)$ and $\Dc_{\rm poly}^{\bu} (M)$ to be 
graded as Lie algebras, i.e. $T_{\rm poly}^i (M) = \{ 
(i+1)$-polyvector fields$\}$ and $\Dc_{\rm poly}^i (M) \sbs 
{\rm Hom}_{\Cb} \, (A^{\ot (i+1)} , A)$.

The formality theorem of Maxim Kontsevich \cite{K} states that 
$T_{\rm poly}^{\bu} (M)$ and $\Dc_{\rm 
poly}^{\bu} (M)$ are quasi-isomorphic as 
differential
graded ($dg$) 
Lie algebras, 
i.e. there exists a $dg$ Lie algebra $\Mc^{\bu}$ and the 
diagram
$$
\matrix{
&&\Mc^{\bu} \cr \cr
&^{\vp_1} \nearrow &&\nwarrow ^{\vp_2} \cr \cr
&T_{\rm poly}^{\bu} &&\Dc_{\rm poly}^{\bu} \cr
}
$$
where $\vp_1$ and $\vp_2$ are maps of the Lie algebras and 
quasi-isomorphisms of the complexes.

In the case $M = \Rb^d$ an explicit 
$L_{\ify}$-quasi-isomorphism $\Uc : T_{\rm poly}^{\bu} 
(\Rb^d) \ra \Dc_{\rm poly}^{\bu} (\Rb^d)$ was constructed.

The cyclic formality conjecture relates to the formality 
theorem like the cyclic complex of an associative algebra $A$ 
relates to its Hochschild complex. It turns out that the 
definition of the {\it cohomological} cyclic complex depends 
on an additional datum -- a trace ${\rm Tr} : A \ra \Cb$ on 
the algebra $A$. In the case $A = C^{\ify} (M)$ it depends on 
a volume form $\Om$ on the manifold $M$.

Let us suppose that the volume form $\Om$ is fixed.

\medskip

\noindent {\bf Definition-lemma.} (cyclic shift operator) 

\noindent For any polydifferential Hochschild cochain $\psi : 
C^{\ify} (M)^{\ot k} \ra C^{\ify} (M)$ there exists a 
polydifferential Hochschild cochain $C(\psi) : C^{\ify} 
(M)^{\ot k} \ra C^{\ify} (M)$ such that for any $k+1$ 
functions $f_1 , \ldots , f_{k+1}$ on $M$ with compact 
support one has:
\begin{equation}
\int_{\Rb^d} \psi (f_1 \ot \ldots \ot f_k) \cdot f_{k+1} 
\cdot \Om = (-1)^k \int_{\Rb^d} C(\psi) (f_2 \ot \ldots \ot 
f_{k+1}) \cdot f_1 \cdot \Om \, .
\label{eq1}
\end{equation}

\noindent {\bf Definition.} (cohomological cyclic complex)
$$
[\Dc_{\rm poly}^{\bu} (M)]_{\rm cycl} = \{ \psi \in \Dc_{\rm 
poly}^{\bu} (M) \mid C(\psi) = \psi \} \, .
$$

\noindent {\bf Lemma.} {\it $[\Dc_{\rm poly}^{\bu} (M)]_{\rm 
cycl}$ is closed under the Hochschild differential and the 
Gerstenhaber bracket.}

\medskip

\noindent See \cite{Sh}, Section~1.3.2 for a proof.

\medskip

We have defined a cyclic analog of the $dg$ Lie algebra $\Dc_{\rm 
poly}^{\bu} (M)$. The cyclic analog of $T_{\rm poly}^{\bu} 
(M)$ is defined as follows: it is $T_{\rm poly}^{\bu} (M) \ot 
\Cb \, [u]$, $\deg u = 2$ with the $\Cb \, [u]$-linear bracket and the 
differential $d_{\rm div} (\g \ot  u^k) = ({\rm div} \, 
\g) \ot u^{k+1}$. The divergence operator ${\rm div} : 
T_{\rm poly}^{\bu} (M) \ra T_{\rm poly}^{\bu - 1} (M)$ is 
defined from the volume form $\Om$ and the de Rham operator:
$$
{\rm div} : T_{\rm poly}^k \build \longra_{}^{\Om} 
\Om^{d-k-1} (M) \build \longra_{}^{d_{\rm DR}} \Om^{d-k} (M) 
\build \longra_{}^{\Om} T_{\rm poly}^{k-1} (M)
$$
(here $d=\dim M$).

We have to prove that $\{T_{\rm poly}^{\bu} \ot \Cb \, [u], 
d_{\rm div} \}$ is actually a $dg$ Lie algebra, i.e. to prove 
that
\begin{equation}
{\rm div} \, [\g_1 , \g_2] = [{\rm div} \, \g_1 , \g_2] \pm 
[\g_1 , {\rm div} \, \g_2] \, .
\label{eq2}
\end{equation}
It follows from the fact that for any volume form $\Om$
\begin{equation}
\pm ( \, {\rm div} \, (\g_1 \wdg \g_2) - ({\rm div} \, \g_1) 
\wdg \g_2 \pm \g_1 \wdg {\rm div} \, \g_2) = [\g_1 , \g_2] \, 
.
\label{eq3}
\end{equation}
Formula (\ref{eq2}) can be obtained from (\ref{eq3}) by the
application of div to both sides and from the identity ${\rm div}^2 
= 0$.

The cyclic formality is the following

\medskip

\noindent {\bf Conjecture.}  {\it The $dg$ Lie algebras
$$
\{ T_{\rm poly}^{\bu} \ot \Cb \, [u] , d_{\rm div} \} \quad 
\hbox{and} \quad \{ [\Dc_{\rm poly}^{\bu} (M)]_{\rm cycl} , 
d_{\rm Hoch} \}
$$
are quasi-isomorphic for any manifold $M$ and volume form 
$\Om$.}

\medskip

This conjecture is due to M.~Kontsevich (unpublished). In 
\cite{Sh} there was constructed (conjecturally) an explicit 
$L_{\ify}$-quasi-isomorphism in the case $M = \Rb^d$.

One can consider also $[T_{\rm poly}^{\bu} (M)]_{\rm div} = \{ 
\g \in T_{\rm poly}^{\bu} (M) \mid {\rm div} \, \g = 0 \}$ as a 
$dg$ Lie algebra with zero differential. The main result of the 
present paper is the following.

\medskip

\noindent {\bf Theorem.}  {\it Let $\Om$ be a constant volume form on $\Rb^d$.
Then the restriction of Kontsevich's 
$L_{\ify}$-quasi-isomorphism $\Uc : T_{\rm poly}^{\bu} (\Rb^d)$ 
$\ra \Dc_{\rm poly}^{\bu} (\Rb^d)$, constructed from the angle 
function
$$
\varphi^h (z,w) = \frac{1}{2i} \, {\rm log} \left( \frac{(z-w) 
(z-\ov w)}{(\ov z - w) (\ov z - \ov w)} \right) \, ,
$$
to the Lie subalgebra $[T_{\rm poly}^{\bu} (\Rb^d)]_{\rm div}$, 
defines an $L_{\ify}$-map $\Uc : [T_{\rm poly}^{\bu} 
(\Rb^d)]_{\rm div} \ra$ \break $[\Dc_{\rm poly}^{\bu} 
(\Rb^d)]_{\rm cycl}$.}

\medskip

In other words, the components $\Uc_k (\g_1 \wdg \cdots 
\wdg \g_k)$ of the Kontsevich $L_{\ify}$-morphism are cyclic, 
if
$$
{\rm div} \, \g_1 = \cdots = {\rm div} \, \g_k = 0 \, 
$$
(with respect to a constant volume form).

\medskip

\noindent {\bf Corollary.}  {\it 
For a Poisson bivector field 
$\pi$ and a constant volume form $\Om$ on $\Rb^d$ such that ${\rm div}_\Om 
\, \pi = 0$, the Kontsevich star-product, constructed from 
$\pi$, is cyclic, i.e.
$$
\int_{\Rb^d} (f * g) \cdot h \cdot \Om = \int_{\Rb^d} (g * h) 
\cdot f \cdot \Om
$$
for any three functions $f,g,h \in C^{\ify} (\Rb^d)$ with 
compact support.}

\medskip

\noindent {\bf Remark.}   Let us note that the complexes 
$[T_{\rm poly}^{\bu} (M)]_{\rm div}$ and $[\Dc_{\rm poly}^{\bu} 
(M)]_{\rm cycl}$ have different cohomology (see \cite{Sh}, 
Section~2), in particular, the $L_{\ify}$-morphism
$$
\Uc : [T_{\rm poly}^{\bu} (\Rb^d)]_{\rm div} \ra [\Dc_{\rm 
poly}^{\bu} (\Rb^d)]_{\rm cycl}
$$
is {\rm not} a quasi-isomorphism.

\medskip

We prove the theorem in the next section and globalize it to the
case of an arbitrary manifold $M$ with a volume form in
Section~\ref{sec3}. 
In particular, we prove that the statement of the corollary
holds for any volume form on $\Rb^d$ but with a star-product
which does not in general coincide with Kontsevich's one.


\section{Geometry of the cyclic formality 
conjecture}\label{sec2}

\subsection{}\label{ssec2.1} In this section we recall 
Kontsevich's construction \cite{K} of $L_{\ify}$-morphism of 
formality, but in a slightly different form. In fact, we 
replace the 2-dimensional group $G^{(2)} = \{ z \mpo az + b \, ; \ 
a,b \in \Rb \, , \ a > 0 \}$ by
 the whole group ${\rm PSL}_2 
(\Rb) = \left\{ z \mpo \frac{az + b}{cz + d} \, ; \ a,b,c,d \in 
\Rb \, , \ \det \left( \matrix{a &c \cr b &d \cr } \right) = 1 
\right\}$. Actually, the group $G^{(2)}$ is exactly the 
subgroup in ${\rm PSL}_2 (\Rb)$, preserving the point $\ify$.

We consider the disk $D^2 = \{ z \mid \vert z \vert \leq 1 \}$ 
instead of the upper half-plane $\Hc$ in \cite{K},
and we
identify it with 
$\Hc\cup\{\infty\}$ by stereographic projection. In
particular the 
group ${\rm PSL}_2 (\Rb)$ of holomorphic transformations of 
$\Cb P^1$ preserving the upper half-plane acts on $D^2$.
 Now an {\it 
admissible graph} is the same as in \cite{K}, but the vertices 
$\{ 1,2, \ldots , n \}$ of the first type are placed in the 
interior of the disk $D^2$, the vertices $\{ \ov 1 , \ov 2 , \ldots 
, \ov m \}$ of the second type are placed on the boundary $S^1 
= \partial D^2$, and $2n + m - 3 \geq 0$. In particular, there 
are no simple loops, and for every vertex $k \in \{ 1,2,\ldots 
, n \}$ of the first type, the set of edges ${\rm Star} \, (k) 
= \{ (v_1 , v_2) \in E_{\G} \mid v_1 = k \}$, starting from $k$, 
is labelled by symbols $(\ell_k^1 , \ldots , \ell_k^{\# {\rm 
Star} (k)})$.

To each such  graph $\G$ with $2n+m-3+\ell$ edges, we attach a 
linear map
$$
\wt{\Uc}_{\G} : \ot^n \, T_{\rm poly}^{\bu} (\Rb^d) \ra 
\Dc_{\rm poly}^{\bu} (\Rb^d) \, [2+\ell - n] \, ,
$$
exactly as in \cite{K}, Section~6.3.

Now we are going to define the weight $W_{\G}$ for a graph $\G$ 
with $n$ vertices of the first type, $m$ vertices of the second 
type, and $2n + m - 3$ edges.

We consider configuration spaces $\wt{\rm Conf}_{n,m} = \{ (p_1 
, \ldots , p_n \, ; \ q_1 , \ldots , q_m) \mid p_i \in {\rm 
Int} \, D^2 \, , \ q_j \in S^1 = \partial D^2 \, , \ p_{i_1} 
\ne p_{i_2}$ for $i_1 \ne i_2$ and $q_{j_1} \ne q_{j_2}$ for 
$j_1 \ne j_2 \}$. We suppose $2n + m \geq 3$. 
Then the group 
${\rm PSL}_2 (\Rb)$
acts freely on $\wt{\rm Conf}_{n,m}$. We set 
$D_{n,m} = \wt{\rm Conf}_{n,m} / {\rm PSL}_2 (\Rb)$; $\dim 
D_{n,m} = 2n + m - 3$. One can construct a compactification 
$\ov{D}_{n,m}$ of the space $D_{n,m}$ exactly as in \cite{K}, 
Section~5, but we will not use it.

Now to each graph $\G$ as above we attach a differential form of 
top degree on the space $D_{n,m}$. Let $\a_1 , \ldots , \a_m$ 
be (positive) real numbers, and let $p,q \in D^2$.

\medskip

\noindent {\bf Definition.} 
Let $\xi_1,\dots,\xi_m\in S^1$.
Let $\vp_i (p,q)$, $i = 1, 
\ldots , m$, be the angle between the geodesic in the 
Poincar\'e metric on $D^2$, connecting $p$ and $q$, and the 
geodesic, connecting $p$ and a point $\xi_i$ on $S^1$ (where 
``sits'' the $i$-th vertex of the second type), see 
Figure~1. 
It is defined modulo $2\pi$.
We set $\vp_{\a_1 , \ldots , \a_m} (p,q) = 
\build \sum_{k=1}^{m} \a_k \, \vp_k (p,q)$. The function 
$\vp_{\a_1 , \ldots , \a_m} (p,q)$ 
depends on the points $p,q \in D^2$ and on the points $\xi_1 , 
\ldots , \xi_m \in S^1$. The differential $d\vp_{\a_1,\ldots,\a_m}
 (p,q)$ is 
a
well-defined 1-form on $\wt{\rm Conf}_{2,m}$. It is ${\rm PSL}_2 
(\Rb)$-invariant, where ${\rm PSL}_2 (\Rb)$ acts simultaneously 
on $p,q$ and on $\xi_1 , \ldots , \xi_m$.

\sevafig{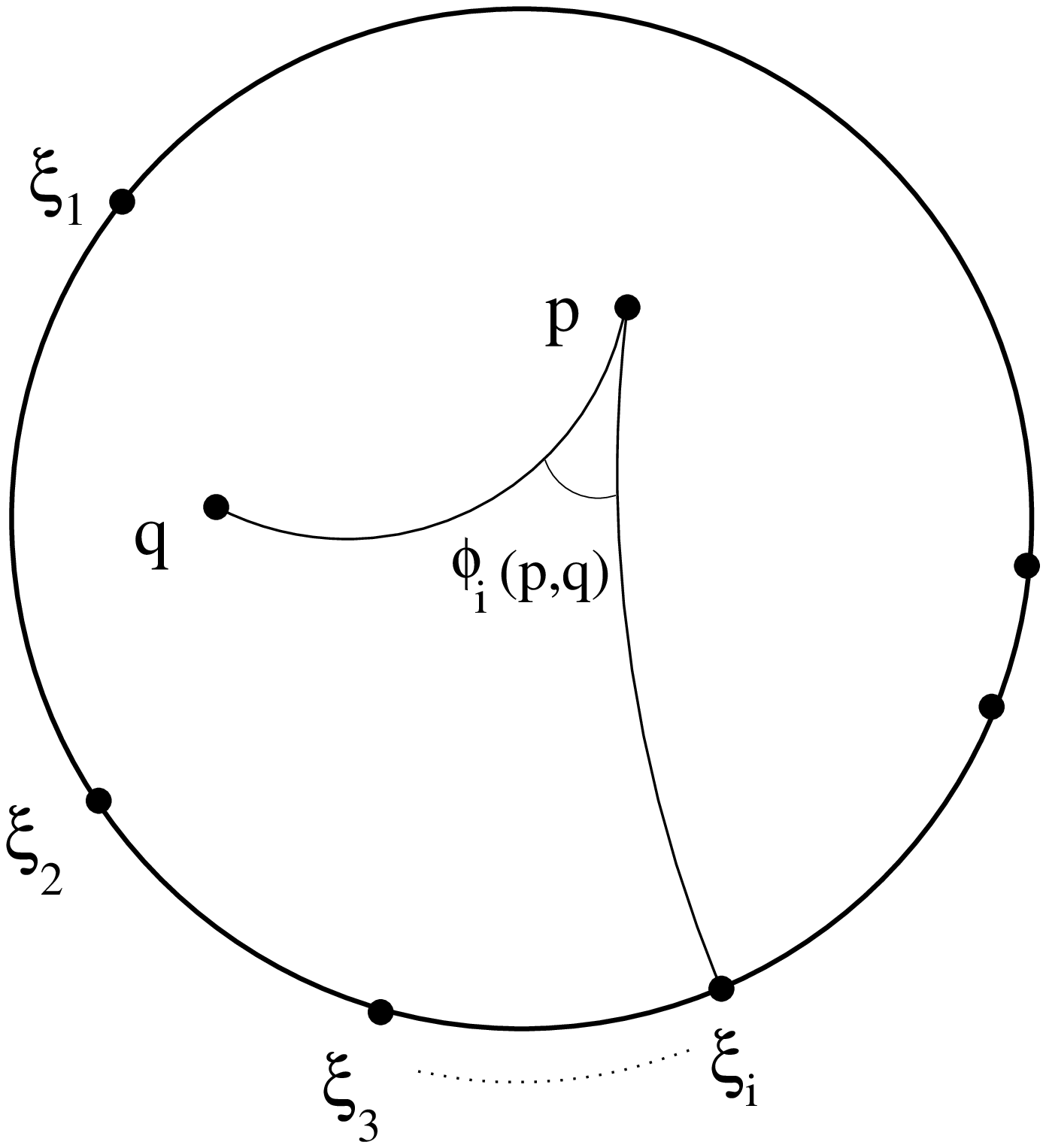}{50mm}

\bigskip

\centerline{Figure 1}

\bigskip

\noindent {\bf Example.} 
If $(\alpha_1,\dots,\alpha_m)=(0,\dots,0,1)$ and $\xi_m=\infty$,
then
$$
\vp_{\a_1,\dots,\a_m} (p,q) = \frac{1}{2i} \, {\rm log} \left( \frac{(p-q) (p - 
\ov q)}{(\ov p - q) (\ov p - \ov q)} \right) \, ,
$$
as in \cite{K}.

\medskip

Now, if the graph $\G$ is ``placed'' on the disk $D^2$, one 
associates to each edge $e$  the 1-form $d \vp_{e; \a_1 , 
\ldots , \a_m}$ on $\wt{\rm Conf}_{n,m}$. 
It is the pull-back of $d \vp_{\a_1 , 
\ldots , \a_m}$ by the map sending the two points 
connected by the edge to $p$,and $q$, and 
$q_j$ to $\xi_j$, $j=1,\dots,m$.
This 1-form induces a 
1-form on the space $D_{n,m}$.

\medskip

\noindent {\bf Definition.} (weight $W_{\G}$)
$$
W_{\G}^{(\a_1 , \ldots , \a_m)} = \prod_{k=1}^{n} \frac{1}{(\# 
\, {\rm Star} \, (k))!} \cdot \frac{1}{(2\pi)^{2n+m-2}} \cdot 
\int_{D_{n,m}^+} \bigwedge_{e \in E_{\G}} d \, \vp_{e;\a_1 , 
\ldots , \a_m} \, .
$$
Here $D_{n,m}^+$ is the connected component consisting of 
configurations for which the points $q_j$, $j = 1,\ldots ,m$, are 
cyclically ordered counterclockwise.

To each map $\Uc_{\G} : \ot^n \, T_{\rm poly}^{\bu} \ra 
\Dc_{\rm poly}^{\bu} [2-n]$ one associates the corresponding 
skew-symmetric map $\wt{\Uc}_T : \bigwedge^n \, T_{\rm 
poly}^{\bu} \ra \Dc_{\rm poly}^{\bu} [2-n]$. We define
$$
\wt{\Uc}_n^{\a_1 , \ldots ,\a_m} = \sum_{m \geq 0} \ \sum_{\G \in 
G_{n,m}} W_{\G}^{\a_1 , \ldots ,\a_m} \ts \wt{\Uc}_{\G}
$$
where $G_{n,m}$ is the set of admissible graphs with $n$ 
vertices of the first type, $m$ vertices of the second type, 
and $2n + m -3$ edges.

\medskip

\noindent {\bf Lemma.} {\it In the case $(\a_1 , \ldots ,\a_m) = 
(0,0,\ldots ,0,1)$
$$
\wt{\Uc}_n^{\a_1 , \ldots ,\a_m} (f_1 \ot \cdots \ot f_m) = 
\Uc_n (f_1 \ot \cdots \ot f_{m-1}) \cdot f_m \, ,
$$
where $\Uc_n$ is the Taylor component of Kontsevich's 
$L_{\ify}$-map (\cite{K}, Section~6).}

\medskip

\noindent {\it Proof.} Using ${\rm PSL}_2 (\Rb)$-invariance, one 
can assume that $f_m$ ``sits'' in the point $\{ \ify \}$, and 
then
the configuration is $G^{(2)}$-invariant. By  
definition of $\vp_{e;\a_1 , \ldots ,\a_m}$, if the graph $\G$ 
contains an edge ending at $\{ \ify \}$, the weight 
$W_{\G}^{\a_1 , \ldots ,\a_m}$
vanishes
 (the corresponding angle 
between two geodesics is zero).~$\square$

\subsection{}\label{ssec2.2}

\noindent {\bf Theorem.} {\it Let $\Om$ be a constant volume form on 
$\Rb^d$, and let $\g_1 , \ldots , \g_n \in [T_{\rm poly}^{\bu} 
(\Rb^d)]_{\rm div}$. Then the integral $\int_{\Rb^d} \, 
\wt{\Uc}_n^{\a_1 , \ldots ,\a_m} (\g_1 \wdg \cdots \wdg \g_n) 
(f_1 \ot \cdots \ot f_m) \cdot \Om$ depends only on the sum 
$\a_1 + \cdots + \a_m$.}

\medskip

The theorem in Section~\ref{sec1} follows from this theorem and 
Lemma~\ref{ssec2.1}, when we set $(\a_1 , \ldots ,\a_m) = (0,0, 
\ldots , 0, 1)$ and $(\a'_1 , \ldots ,\a'_m) = (1,0, \ldots 
,0)$.

\medskip

\noindent {\it 
Proof of the theorem.} We consider the case $(\a_1 , 
\ldots ,\a_m) = (0,0, \ldots , 0, 1)$ and $(\a'_1 , \ldots 
,\a'_m) = (1,0, \ldots ,0)$; the general case is analogous.

\medskip

\noindent {\bf Key-lemma.}
{\it
$$
d \, \vp_{\a_1 , \ldots ,\a_m} (p,q) = d \, \vp_{\a'_1 , \ldots 
,\a'_m} (p,q) + 
d\Phi_{\{ \a \} , \{ \a' \}} (p)
$$
where the 1-form 
$d\Phi$ does not depend on $q$.}

\sevafig{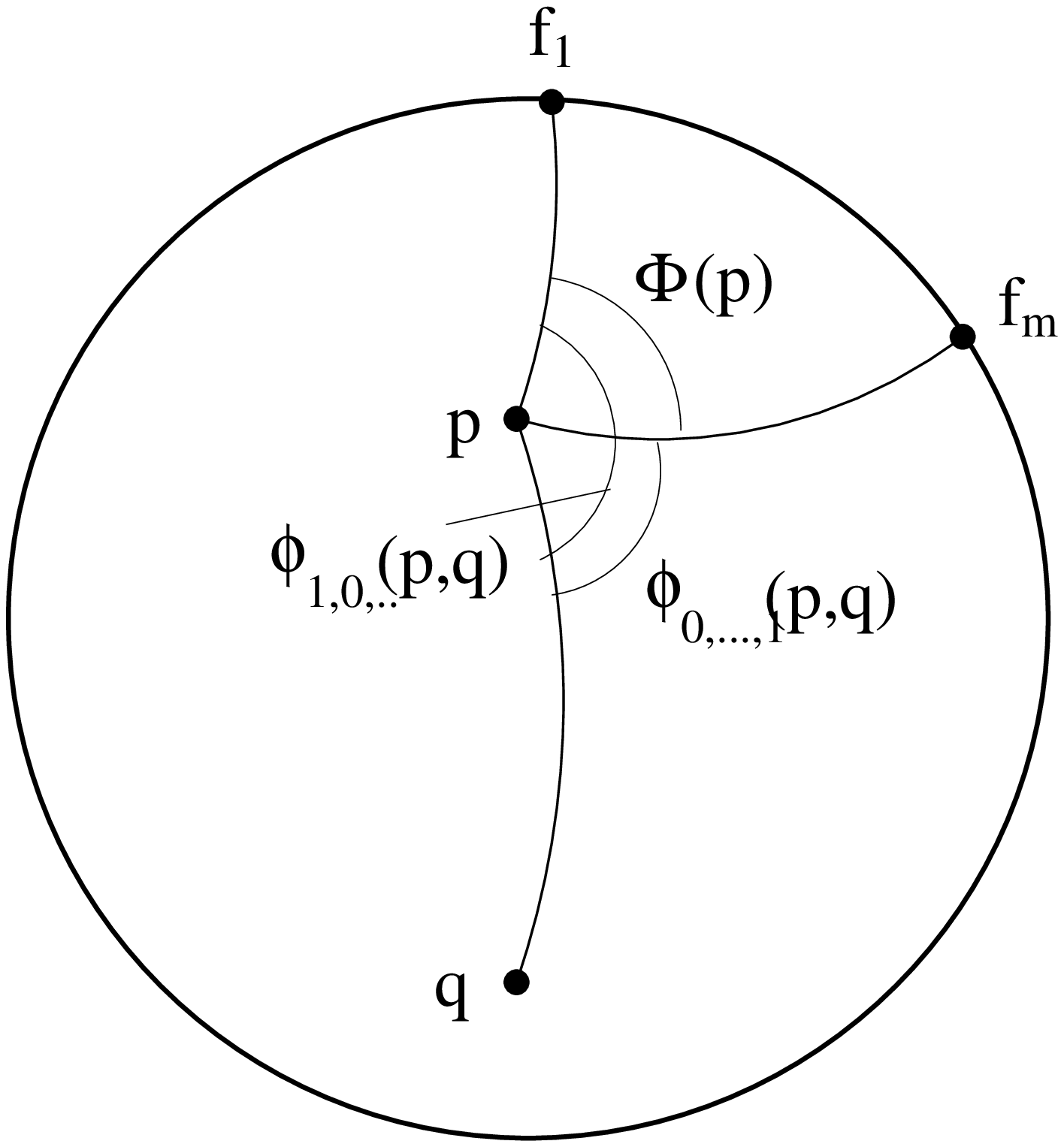}{50mm}
\medskip
\centerline{Figure 2}
\noindent {\it Proof.} It follows directly from the additivity 
of the angle function, see Figure~2.~$\square$

\medskip

We proceed to prove Theorem~\ref{ssec2.2}. The weight 
$W_{\G}^{\a_1 , \ldots ,\a_m}$ is the integral of the wedge 
product $\build\bigwedge_{e \in E_{\G}}^{} d \vp_{e;\a_1 , 
\ldots ,\a_m}$. Now $d \vp_{e;\a'_1 , \ldots ,\a'_m} = d 
\vp_{e;\a_1 , \ldots ,\a_m} + 
d\Phi_e (p)$. Then, by the 
skew-symmetry,
\begin{eqnarray}
&&\bigwedge_{e \in {\rm Star} (k)} d \vp_{e;\a'_1 , \ldots 
,\a'_m} = \\
&&\sum_{e \in {\rm Star} (k)} \left( \pm \bigwedge_{\ov e \in 
{\rm Star} (k) \backslash e} d \vp_{\ov e;\a_1 , \ldots ,\a_m} 
\right) \wdg 
d\Phi_e 
+ \bigwedge_{e \in {\rm Star} (k)} d 
\vp_{e;\a_1 , \ldots ,\a_m} \, . \label{eq4} \nonumber
\end{eqnarray}
Let us denote by $\om_e$ the form
$$
\left( \bigwedge_{\ov e \in {\rm Star} (k) \backslash e} d 
\vp_{\ov e ; \a_1 , \ldots ,\a_m} \right) \wdg 
d\Phi_e
\wdg 
\bigwedge_{\wt e \notin {\rm Star} (k)} d \vp_{\wt e ; \a_1 , 
\ldots ,\a_m}
$$
($e \in {\rm Star} (k)$).

If $\G' = (\G \backslash \{ e \} ) \coprod \{ e' \}$ where the 
edges $e$ and $e'$ have the same start-point, then
\begin{equation}
\int_{D_{n,m}} \om_e^{\G} = \int_{D_{n,m}} \om_{e'}^{\G'} \, . 
\label{eq5}
\end{equation}

Now it is sufficiently to prove the following

\medskip

\noindent {\bf Lemma.} {\it Let $\ov{\G}$ be a graph with $n$ 
vertices of the first type, $m$ vertices of the second type, and 
$2n+m-4$ edges; and let $\g_1 , \ldots , \g_n$ be arbitrary 
polyvector fields. Let $G_{\ov \G}^i$ be the set of all the 
graphs, obtained from the graph $\ov{\G}$ by addition of any 
edge $e$ such that the start-point of $e$ is $\{ i \}$. Then}
\begin{eqnarray}
&&\int_{\Rb^d} \ \sum_{\G \in G_{\ov{\G}}^i} \Uc_{\G} (\g_1 \wdg 
\cdots \wdg \g_n) (f_1 \ot \cdots \ot f_m) \cdot \Om \\
& =& \pm \int_{\Rb^d} \Uc_{\ov \G} (\g_1 \wdg \cdots \wdg {\rm 
div} \g_i \wdg \cdots \wdg \g_n) (f_1 \ot \cdots \ot f_m) \cdot 
\Om \, . \nonumber \label{eq6}
\end{eqnarray}

This is a standard result. The simplest version of it is the 
following: let $\xi$ be a vector field on $\Rb^d$, then
$$
\int_{\Rb^d} \xi (f) \cdot \Om = - \int {\rm div} \, (\xi) \cdot 
f \cdot \Om \, .
$$
Indeed, let us suppose that $\Om = dx_1 \wdg \cdots \wdg d 
x_d$, then, if $\xi = \sum a_i (x_1 , \ldots , x_d) \, 
\frac{\partial}{\partial x_i}$, one has:
$$
\int_{\Rb^d} \xi (f) = \int \sum a_i \cdot \partial_i (f) = 
\int \sum \partial_i (a_i \cdot f) - \int \sum \partial_i (a_i) 
\cdot f \, .
$$
The first summand is equal to 0, and the second summand is $- 
\int ({\rm div} \, \xi) \cdot f$.

Now Theorem~\ref{ssec2.2} follows from the last lemma and 
formula (\ref{eq5}).

Theorem~\ref{ssec2.2} is proven, and Theorem~\ref{sec1} follows 
now from Theorem~\ref{ssec2.2} and Lemma~\ref{ssec2.1}.~$\square$

\medskip

\noindent {\bf Remark.} For a general volume form the proof fails, because
(6) is not true.

\medskip

\noindent {\bf Remark.} 
We have proved the theorem for 
very special choice of angle function (in the sense of 
\cite{K}, Section~6.2), namely, we use the harmonic angle 
function, which also appears in the QFT approach to the 
formality theorem \cite{CF}. It seems that Theorem~\ref{sec1} 
is not true for any other choice.
\medskip

\newcommand{\lp}{{\scriptstyle(}}
\newcommand{\rp}{{\scriptstyle)}}

\noindent{\bf Remark.}
The cyclicity of the star product may be heuristically 
understood, 
for $M=\Rb^d$ with volume form $\Omega=dx_1\wdg\cdots\wdg dx_d$,
in ``physical'' terms as follows: in \cite{CF} the
Kontsevich star
product was described as the Feynman perturbation
expansion of a path integral formula
$f*g(x)=\int_{X\lp p_3\rp =x}\exp(\frac i\hbar S(\hat X))
f(X\lp p_1\rp )g(X\lp p_2\rp )d\hat X$, for a certain action functional
$S$ on the space of bundle homomorphisms $\hat X:TD\to T^*M$,
with base map $X:D\to M$,
from the tangent bundle of an oriented disk $D$ to the cotangent
bundle of $M$. The points $p_1,p_2,p_3$ are any three cyclically
ordered points on the boundary of $D$. One may more
generally consider for any three functions $f,g,h\in C^\infty(M)$
the correlation function
\[
\langle f,g,h\rangle=\int e^{\frac i\hbar S(\hat X)}
f(X\lp p_1\rp )g(X\lp p_2\rp )h(X\lp p_3\rp )d\hat X
\]
which looks invariant under cyclic permutations of $f,g,h$, since
the action is invariant under orientation preserving diffeomorphisms
of $D$. Moreover, if $M=\Rb^d$, we may write this integral as the integral over
the maps with $x=\sum\alpha_i X\lp p_i\rp$ fixed, and then over the
position of the ``center of mass'' $x$.
Naively, this is independent of $\alpha_i$ with
$\sum\alpha_i=1$. In particular, with $\alpha=(0,0,1)$ and $(1,0,0)$,
we obtain
\[
\int_{\Rb^d} f* g\cdot h\, dx_1\wdg\cdots\wdg dx_d=\int_{\Rb^d}
 g*h\cdot f\, dx_1\wdg\cdots\wdg dx_d.
\]
However there is an anomaly, meaning that the symmetry under diffeomorphisms
of the disk is not a symmetry of the path integral. Technically
this follows from the fact that the
regularization of amplitudes of Feynman diagrams involving
tadpoles (edges with both ends at the same vertex) cannot be chosen
in an invariant way, see \cite{CF}. But the tadpoles correspond
to bidifferential operators involving the divergence of the Poisson 
bivector field.
Therefore the anomalous terms vanish 
for divergence free Poisson bivector fields and the above argument applies.

\medskip

\section{Globalization}\label{sec3}

\subsection{}\label{ssec3.1} Here we prove the following

\medskip

\noindent {\bf Theorem.} {\it For any smooth manifold $M$ and any 
volume form $\Om$ on $M$, there exists an $L_{\ify}$-morphism
$$
\Uc : [T_{\rm poly}^{\bu} (M)]_{\rm div} \ra [\Dc_{\rm 
poly}^{\bu} (M)]_{\rm cycl} \, ,
$$
such that its first Taylor component, $\Uc_1$, coincides with 
the Hochschild-Kostant-Rosenberg map.}

\medskip

The proof follows basically the same line as the proof of 
the globalization of the formality theorem in \cite{K}, 
Section~7. Let us note that in our case the $L_{\ify}$-morphism 
$\Uc : [T_{\rm poly}^{\bu} (M)]_{\rm div} \ra [\Dc_{\rm 
poly}^{\bu} (M)]_{\rm cycl}$ is {\it not} an 
$L_{\ify}$-quasi-isomorphism. 

For any $d$-dimensional manifold $M$, there exists an 
infinite-dimensional manifold $M^{\rm coor}$, the manifold of 
formal coordinate systems on $M$. The main property of the 
manifold $M^{\rm coor}$ is that there exists a $W_d$-valued 
1-form $\om$ on $M^{\rm coor}$, satisfying the Maurer-Cartan 
equation $d\om + \frac{1}{2} \, [\om , \om] = 0$ (here $W_d = 
{\rm Vect} \, (\Rb_{\rm formal}^d)$ is the Lie algebra of formal 
vector fields on $\Rb^d$). In other notations, there exists a 
map $T [1] \, M^{\rm coor} \ra W_d [1]$, which is a map of 
$Q$-manifolds (we refer to \cite{K} for basic definitions on 
$Q$-manifolds). We need to modify this construction in the case 
when the manifold $M$ is equipped with a volume form $\Om$.

We define an infinite-dimensional manifold $M_{\Om}^{\rm coor}$ 
for any $d$-dimensional manifold $M$ with volume form $\Om$, as 
follows. A point of the manifold $M_{\Om}^{\rm coor}$ is a map 
$\vp : \Rb_{\rm formal}^d \ra M$, $\vp (0) = x \in M$ such that 
$\vp^* \, \Om = dx_1 \wdg \cdots \wdg d x_d$. Let us denote by 
$[W_d]_{\rm div}$ the Lie algebra of formal vector fields on 
$\Rb^d$ with zero divergence (with respect to the standard 
volume form $dx_1 \wdg \cdots \wdg d x_d$ on $\Rb^d$). More 
explicitly, $[W_d]_{\rm div} = \left\{ \sum a_i (x_1, \ldots ,x_n) \, \partial_i 
\mid \sum \partial_i \, a_i = 0 \right\}$ . 

\medskip

\noindent {\bf Remark.} The volume form $\Om$ on $M$ is constant in the coordinates
in a neighbourhood of the point $x$ obtained through the map $\vp$ from the affine coordinates on 
$\Rb_{\rm formal}^d$.

\medskip

There exists a $[W_d]_{\rm 
div}$-valued 1-form $\om_{\rm div}$ on $M_{\Om}^{\rm coor}$, which satisfies the 
Maurer-Cartan 
equation $d \, \om_{\rm div} + \frac{1}{2} \, [\om_{\rm div}, 
\om_{\rm div}] =0$. In the other terms, there exists a map of 
$Q$-manifolds $T [1] \, M_{\Om}^{\rm coor} \ra [W_d]_{\rm div} 
\, [1]$. The Lie algebra $sl_d = \left\{ \sum a_{ij} \, x_i \, 
\partial_j \mid a_{ij} \in \Cb \, , \ tr (a_{ij}) = 0 \right\} 
\sbs [W_d]_{\rm div}$ acts on $M_{\Om}^{\rm coor}$ (as well 
as the whole Lie algebra $[W_d]_{\rm div}$), and this action 
can be integrated to an action of the group ${\rm SL}_d$.

\medskip

\noindent {\bf Lemma.} {\it The fibers of the natural bundle 
$M_{\Om}^{\rm coor} / {\rm SL}_d \ra M$ are contractible.}

\medskip

\noindent {\it Proof.} Let $\vp : \Rb_{\rm formal}^d \ra 
\Rb_{\rm formal}^d$, $\vp (0) = 0$, be a formal diffeomorphism, 
preserving the volume form $\Om = dx_1 \wdg \cdots \wdg dx_d$. 
Then the map $\vp_{\hbar} : \Rb_{\rm formal}^d \ra \Rb_{\rm 
formal}^d$,
$$
\vp_{\hbar} (x_1 , \ldots , x_d) = \frac{\vp (\hbar \, x_1 , 
\ldots , \hbar \, x_d)}{\hbar} \, , \ \hbar \ne 0 \, ,
$$
also preserves the volume form, 
and so does 
its limit when 
$\hbar \ra 0$. It is clear that this limit is a linear map 
$\vp_0 : \Rb_{\rm formal}^d \ra \Rb_{\rm formal}^d$, and, 
therefore, $\vp_0 \in {\rm SL}_d$. We have constructed a 
retraction of a fiber of the bundle $M_{\Om}^{\rm coor} \ra M$ 
on the space ${\rm SL}_d$.~$\square$

\bigskip

We  use this lemma and the following 
properties of the $L_{\ify}$-morphism $\Uc : [T_{\rm 
poly}^{\bu} (\Rb^d)]_{\rm div} \ra [\Dc_{\rm poly}^{\bu} 
(\Rb^d)]_{\rm cycl}$:
\begin{itemize}
\item[P1)] $\Uc$ can be defined for $\Rb_{\rm formal}^d$ as 
well;
\item[P2)] for any $\xi \in [W_d]_{\rm div}$ we have
$$
\Uc_1 (m_T (\xi)) = m_D (\xi)
$$
(here $m_T : [W_d]_{\rm div} \ra [T_{\rm poly}^{\bu} (\Rb_{\rm 
formal}^d)]_{\rm div}$ and $m_D : [W_d]_{\rm div} \ra [\Dc_{\rm 
poly}^{\bu} (\Rb_{\rm formal}^d)]_{\rm cycl}$ are the canonical 
maps);
\item[P3)] $\Uc$ is ${\rm SL}_d$-equivariant;
\item[P4)] for any $k \geq 2$, $\xi_1 , \ldots ,\xi_k \in 
[W_d]_{\rm div}$ one has
$$
\Uc_k (m_T (\xi_1) \ot \cdots \ot m_T (\xi_k)) = 0 \, ;
$$
\item[P5)] for any $k \geq 2$, $\xi \in sl_d \sbs [W_d]_{\rm 
div}$ and for any $\eta_2 , \ldots ,\eta_k \in T_{\rm 
poly}^{\bu} (\Rb_{\rm formal}^d)$ one has:
$$
\Uc_k (m_T (\xi) \ot \eta_2 \ot \cdots \ot \eta_k) = 0 \, .
$$
\end{itemize}

The properties P1)--P5) are cited from \cite{K}, Section~7, 
where they were proven for the (a bit 
stronger) case of the 
$L_{\ify}$-map of formality $\Uc : T_{\rm poly}^{\bu} (\Rb^d) 
\ra \Dc_{\rm poly}^{\bu} (\Rb^d)$,  we just replace the Lie 
algebra $W_d$ by $[W_d]_{\rm div}$ and its subalgebra $gl_d 
\sbs W_d$ by the subalgebra $sl_d \sbs [W_d]_{\rm div}$. 
According to Remark 3.1,
for the globalization it is sufficient to know the local result only 
for a constant volume form $\Om$.

The theorem can be deduced from the properties P1)-P5) 
exactly in the same way as it is done in the case of 
Formality in \cite{K}, Section~7.

\subsection{Consequences}\label{ssec3.2}

Here we consider some consequences of Theorem 3.1.

\medskip

\noindent {\bf Corollary 1.} {\it Let $M$ be a Poisson manifold 
(with the bivector field $\pi$, and let $\Om$ be any volume 
form $\Om$ on $M$ such that ${\rm div}_{\Om} \, \pi = 0$.
Then there exists a star-product on $C^{\ify} (M)$ such that 
for any three functions $f,g,h$ with compact support one has:
$$
\int_M (f * g) \cdot h \cdot \Om = \int_M (g * h) \cdot f \cdot 
\Om \, .
$$}

\noindent {\it Proof.} We apply to the $L_{\ify}$-morphism $\Uc 
: [T_{\rm poly}^{\bu} (M)]_{\rm div} \ra [\Dc_{\rm poly}^{\bu} 
(M)]_{\rm cycl}$, constructed above, the following general 
statement. For an $L_{\ify}$-morphism $\Fc : \Gc_1^{\bu} \ra 
\Gc_2^{\bu}$ between two $dg$ Lie algebras and a solution $\pi$ 
of the Maurer-Cartan equation in $(\Gc_1^{\bu})^1$, the formula
$$
\Fc (\pi) = \Fc_1 (\pi) + \frac{1}{2} \, \Fc_2 (\pi , \pi) + 
\frac{1}{6} \, \Fc_3 (\pi , \pi , \pi) + \cdots
$$
defines a solution of the Maurer-Cartan equation in 
$(\Gc_2^{\bu})^1$. We apply this construction to the solution 
$\hbar \pi$ of the Maurer-Cartan equation in $[T_{\rm 
poly}^{\bu} (M)]_{\rm div}$.~$\square$

\bigskip

To deduce the Connes-Flato-Sternheimer conjecture from 
Corollary~1, one needs to prove that $1 * f = f * 1 = f$ for any 
function $f$. Locally it is true (for the $L_{\ify}$-morphism 
$\Uc : [T_{\rm poly}^{\bu} (\Rb^d)]_{\rm div} \ra [\Dc_{\rm 
poly}^{\bu} (\Rb^d)]_{\rm cycl}$, constructed in 
Section~\ref{sec2}), therefore it is true also globally. So, we 
have proved

\medskip

\noindent {\bf Corollary 2.} (generalized 
Connes-Flato-Sternheimer conjecture) {\it For the star-product of 
Corollary~1 one has:
$$
\int_M (f * g) \cdot \Om = \int_M f \cdot g \cdot \Om \, .
$$
}

\noindent {\it Proof.} Put $h=1$ and use that $g * 1 = 
g$.~$\square$

\bigskip

\noindent {\bf Acknowledgments.} 
The first author thanks Alberto Cattaneo for useful comments
and discussions.
The second author is grateful 
to Maxim Kontsevich and to Andrej Lossev for many discussions.
We thank Jim Stasheff for commenting on the manuscript and numerous
corrections.
We are thankful to Mme 
   C\'ecile Gourgues for the quality typing of this text.

\vfill
\eject

\noindent Giovanni Felder

\noindent Department of Mathematics

\noindent ETH-Zentrum

\noindent 8092 Z\"urich

\noindent SWITZERLAND

\smallskip

\noindent e-mail address: felder@math.ethz.ch

\vglue 1cm

\noindent Boris Shoikhet

\noindent IHES

\noindent 35, route de Chartres

\noindent 91440 Bures-sur-Yvette

\noindent FRANCE

\smallskip

\noindent e-mail address: boris@ihes.fr, borya@mccme.ru

\end{document}